\documentclass [12pt]{article}
\usepackage{graphicx,amssymb,amsfonts,latexsym,amsmath,amsthm,times}
\usepackage{epsfig,fancyhdr,color,ulem}
\usepackage[english]{babel}
\numberwithin{equation}{section}

\newtheorem{theorem}{Theorem}

\newtheorem{example}[theorem]{Example}

\newtheorem{lemma}[theorem]{Lemma}

\newtheorem{proposition}[theorem]{Proposition}
\newtheorem{remark}[theorem]{Remark}

\newcommand{\eqdef}{{\ \stackrel{\mathrm{def}}{=}\ }}
\newcommand{\A}{\mathcal{A}}
\newcommand{\B}{\mathcal{B}}
\newcommand{\C}{\mathcal{C}}
\newcommand{\D}{\mathcal{D}}
\newcommand{\Fe}{\mathcal{F}_\varepsilon}
\newcommand{\ud}{\mathrm{d}}
\newcommand{\ra}{\rightarrow}

\begin{document}

\footnotesize {\flushleft \mbox{\bf \textit{Math. Model. Nat.
Phenom.}}}
 \\
\mbox{\textit{{\bf Vol. X, No. X, 2009, pp. X-XX}}}

\medskip


\thispagestyle{plain}

\vspace*{2cm} \normalsize \centerline{\Large \bf Semigroup analysis of structured parasite populations} 

\vspace*{1cm}

\centerline{\bf J\'{o}zsef Z. Farkas$^1$\footnote{Corresponding
author. E-mail: jzf@maths.stir.ac.uk}, Darren M. Green$^2$ and Peter Hinow$^3$ }

\vspace*{0.5cm}

\centerline{$^1$Department of Computing Science and Mathematics} 
\centerline{University of Stirling, FK9 4LA, Scotland UK}

\centerline{$^2$Institute of Aquaculture, University of Stirling, FK9 4LA, Scotland, UK}

\centerline{$^3$Department of
Mathematical Sciences, University of Wisconsin -- Milwaukee} 
\centerline{P.O.~Box 413, Milwaukee, WI 53201, USA}


\vspace*{1cm}

\noindent {\bf Abstract.}
Motivated by  structured parasite populations in aquaculture we consider a class of size-structured population models, 
where individuals may be recruited into the population with distributed states at birth. The mathematical model which describes the evolution of such a population is a first-order nonlinear 
partial integro-differential equation of hyperbolic type. First, we use positive perturbation arguments and utilise results from the spectral 
theory of semigroups to establish conditions for the existence of a positive equilibrium solution of our model. Then, we formulate conditions that 
guarantee that the linearised system is governed by a positive quasicontraction semigroup on the biologically 
relevant state space. We also show that the governing linear semigroup is eventually compact, hence growth properties of the semigroup are determined by 
the spectrum of its generator. In the case of a separable fertility function, we deduce a characteristic equation, and investigate the stability 
of equilibrium solutions in the general case using positive perturbation arguments.
\vspace*{0.5cm}




\vspace*{1cm}

\setcounter{equation}{0}
\section{Introduction}

In this paper, we study the following partial integro-differential equation
\begin{align}
\frac{\partial}{\partial t}p(s,t)+\frac{\partial}{\partial s}\left(\gamma(s,P(t))p(s,t)\right)&=-\mu(s,P(t))p(s,t)+\int_0^m\beta(s,y,P(t))p(y,t)\,\ud y,\label{equation} \\
\gamma(0,P(t))p(0,t)&=0,\label{boundary} \\
p(s,0)&=p_0(s),\ \ P(t)=\int_0^m p(s,t)\,\ud s\label{initial}.
\end{align}
Here the function $p=p(s,t)$ denotes the density of individuals of size (or other developmental stage) $s$ at time $t$ with $m$ being the finite maximal size any individual may reach in its lifetime. Vital rates $\mu\ge0$ and $\gamma\ge0$ denote the mortality and growth rates of individuals, 
respectively, and both depend on both size $s$ and on the total population size $P(t)$. 
It is assumed that individuals may have different sizes at birth and therefore  $\beta(s,y,\,\cdot\,)$ denotes the rate at which 
individuals of size $y$ give rise to individuals of size $s$. The non-local integral term in \eqref{equation} 
represents reproduction of the population without external driving of the population through immigration.  We make the following general assumptions on the model ingredients
\begin{align}
    & \mu\in C^1([0,m]\times [0,\infty)),\quad \beta\in C^1([0,m]\times [0,m]\times [0,\infty))\nonumber \\
&  \beta,\,\mu\geq 0,\quad \gamma\in C^1([0,m]\times [0,\infty)),\quad\gamma >0.\label{assumptions}
\end{align}

Our motivation to investigate model \eqref{equation}-\eqref{initial} is the modelling of structured parasite populations in aquaculture. 
In particular we are interested in parasites of farmed and wild salmonid fish that have particular relevance 
both industrially and commercially to the UK. These species are subject to parasitism from a number of copepod (crustacean) parasites of the family Caligidae. 
These sea louse parasites are well studied with a large literature: below we draw attention to some recent key review papers.  
Sea lice cause reduced growth and appetite, wounding, and susceptibility to secondary 
infections \cite{costello}, resulting in significant damage to crops and therefore they are economically important. 
For salmon, louse burden in excess of 0.1 lice per gram of fish can be considered pathogenic \cite{costello}. 
The best studied species is \textit{Lepeophtheirus salmonis}, principally a parasite of salmonids and frequent parasite on 
British Atlantic salmon (\textit{Salmo salar}) farms \cite{tully}. It also infects sea trout (\textit{Salmo trutta}) and rainbow trout 
(\textit{Oncorhynchus mykiss}). The life history of the parasite is direct, with no requirement for intermediate hosts. 
It involves a succession of ten distinct developmental stages, separated by moults, from egg to adult. 
Initial \textit{naupliar} and \textit{copepodid} stages are free living and planktonic. Following attachment of the 
infectious copepodid to a host, the parasite passes through four \textit{chalimus} stages that are firmly attached to the host, 
before entering sexually dimorphic \textit{pre-adult} and \textit{adult} stages where the parasite can once again move over the 
host surface and transfer to new hosts.

The state of the art for population-level modelling of \textit{L. salmonis} is represented by Revie \textit{et al}.~\cite{revie}. 
These authors presented a series of delay-differential equations to model different life-history stages and parameterised the model 
using data collected at Scottish salmon farms. A similar compartmental model was proposed by Tucker \textit{et al}. \cite{tucker}. 
The emphasis of these papers was not however, in analytical study, but on numerical simulation and parameterisation using field \cite{revie} 
and laboratory \cite{tucker} data. An earlier model by Heuch \& Mo \cite{heuch} investigated the infectivity, in term of \textit{L. salmonis} egg production, 
posed by the Norwegian salmon industry, using a simple deterministic model. Other authors have considered the potential for long-distance dispersal 
of mobile parasite stages through sea currents \cite{murray}, looking at Loch Shieldaig in NW Scotland, a long-term study site for sea louse research.

In this paper, we focus on the dynamics of individuals at the chalimus to adult
stages. Though individuals pass through a series of discrete growth stages by
moulting, this outward punctuated growth disguises a physiologically more smooth
growth process in terms of the accumulation of energy, and by `size' in this
paper we presume accumulation of energy, rather than physical dimension. Sea
lice reproduce sexually; however at the chalimus stage individuals are not yet
sexually differentiated. Fertility rates thus must be considered as applying to
the population as a whole, rather than as is usually the case the female
fraction of the population.
Individuals entering the first chalimus stage from the non-feeding planktonic stages are distributed over different sizes, 
hence we have the zero influx boundary condition \eqref{boundary} and the recruitment term in \eqref{equation}. 
Our aim here is to present a preliminary step towards the analysis 
of the more complex problem of modelling the whole life cycle of sea lice by giving a mathematical treatment of a quite general scramble competition model with distributed states-at-birth. 
We use the term scramble competition to describe the scenario where individuals have equal chance when competing for resources 
such as food (see e.g. \cite{CUS}). Therefore all vital rates, i.e. growth, fertility and mortality depend on the 
total population size of competitors. In other populations, such as a tree population or a cannibalistic population, 
there may be a natural hierarchy among individuals of different sizes, which results in mathematical models incorporating infinite-dimensional nonlinearities, 
see e.g. \cite{FH,FH2}.  The analysis presented in this paper could be extended to these type of models and also to other models such as those that involve a different type of recruitment term. 

Here, we consider the asymptotic behaviour of solutions of model \eqref{equation}-\eqref{initial}. 
Our analysis is based on linearisation around equilibrium solutions (see e.g.~\cite{FH,PR1}) and utilises well-known results from 
linear operator theory that can be found for example in the excellent books \cite{AGG,CH,NAG}. We also utilise some novel 
ideas on positive perturbations of linear operators. For basic concepts and results from the theory of structured population dynamics 
we refer the interested reader to \cite{CUS,I,MD,WEB}. 

Traditionally, structured population models have been formulated as partial differential equations for population densities. However, 
the recent unified approach of Diekmann \textit{et al.}, making use of the rich theory of delay and integral equations, 
has been resulted in significant advances. The Principle of Linearised Stability has been proven in \cite{D1,D2} 
for a wide class of physiologically structured population models formulated as delay equations (or abstract integral equations). 
It is not clear yet whether the models formulated in \cite{D1,D2} as delay equations are equivalent to those formulated 
as partial differential equations.

In the remarkable paper \cite{CS2}, Calsina and Salda\~{n}a studied the well-posedness of a very general size-structured model with distributed states-at-birth. 
They established the global existence and uniqueness of solutions utilising results from the theory of nonlinear evolution equations. 
Model \eqref{equation}-\eqref{initial} is a special case of the general model treated in \cite{CS2}, however, in \cite{CS2} qualitative 
questions were not addressed.  In contrast to \cite{CS2}, our paper focuses on the existence and local asymptotic stability of equilibrium solutions of 
system \eqref{equation}-\eqref{initial} with particular regards to the effects of distributed states-at-birth compared to more simple 
models we addressed previously, e.g. in \cite{FH}. First, we establish conditions in Theorem \ref{equilibrium2} that 
guarantee the existence of equilibrium solutions, in general. Then, we show  in Theorem \ref{posgen} that a positive quasicontraction semigroup 
describes the evolution of solutions of the system linearised at an equilibrium solution. Next, we establish a further 
regularity property in Theorem \ref{reg} for the governing linear semigroup, which allows one to investigate the stability of 
positive equilibrium solutions of \eqref{equation}-\eqref{initial}. We use rank-one perturbations of the general recruitment term to arrive 
at stability/instability conditions for the equilibria. Finally we briefly discuss the positivity of the governing linear semigroup.


\vspace*{0.5cm}
\setcounter{equation}{0}
\section{Existence of equilibrium solutions}

Model \eqref{equation}-\eqref{initial} admits the trivial solution. If we look for positive time-independent solutions of \eqref{equation}-\eqref{initial} we arrive at the following integro-differential equation
\begin{align}
& \gamma(s,P_*)p_*'(s)+\big(\gamma_s(s,P_*)+\mu(s,P_*)\big)p_*(s)=\int_0^m\beta(s,y,P_*)p_*(y)\,\ud y \label{stateq}\\
& \gamma(0,P_*)p_*(0)=0,\quad P_*=\int_0^m p_*(s)\,\ud s. \label{statbound}
\end{align}

\subsection{Separable fertility function}

In the special case of
\begin{equation}\label{specbeta}
\beta(s,y,P)=\beta_1(s,P)\beta_2(y),\quad s,y\in [0,m],\quad P\in \textcolor{blue}{(}0,\infty),
\end{equation}
where the distribution of offspring sizes is dependent upon the level of competition $P$, but the mature size at which individuals reproduce is not, equation \eqref{stateq} reduces to
\begin{equation}\label{specstateq}
\gamma(s,P_*)p_*'(s)+\big(\gamma_s(s,P_*)+\mu(s,P_*)\big)p_*(s)=\beta_1(s,P_*)\overline{P}_*,
\end{equation}
where
\begin{equation*}
\overline{P}_*=\int_0^m\beta_2(y)p_*(y)\,\ud y.
\end{equation*}
The solution of \eqref{specstateq} satisfying the initial condition in \eqref{statbound} is readily obtained as
\begin{equation}
p_*(s)=\overline{P}_*F(s,P_*)\int_0^s\frac{\beta_1(y,P_*)}{F(y,P_*)\gamma(y,P_*)}\,\ud y,\label{statsol}
\end{equation}
where
\begin{equation*}
F(s,P_*)=\exp\left\{-\int_0^s\frac{\gamma_s(y,P_*)+\mu(y,P_*)}{\gamma(y,P_*)}\,\ud y\right\}.
\end{equation*}
Multiplying equation \eqref{statsol} by $\beta_2$ and integrating from $0$ to $m$ yields the following necessary condition for the existence of a positive equilibrium solution
\begin{equation}
1=\int_0^m \beta_2(s)F(s,P_*)\int_0^s\frac{\beta_1(y,P_*)}{F(y,P_*)\gamma(y,P_*)}\,\ud y\,\ud s.\label{rep}
\end{equation}
Therefore we define a net reproduction function $R$ as follows
\begin{equation}
R(P)=\int_0^m\int_0^s\frac{\beta_1(y,P)\beta_2(s)}{\gamma(s,P)}\exp\left\{-\int_y^s\frac{\mu(z,P)}{\gamma(z,P)}\,\ud z\right\}\,\ud y\,\ud s.\label{netrep}
\end{equation}
It is straightforward to show that for every positive value $P_*$ for which $R(P_*)=1$ holds, formula \eqref{statsol} yields a unique positive stationary solution $p_*$, where $\overline{P}_*$ may be determined from equation \eqref{statsol} as
\begin{equation*}
\overline{P}_*=\frac{P_*}{\int_0^mF(s,P_*)\int_0^s\frac{\beta_1(y,P_*)}{F(y,P_*)}\,\ud y\,\ud s}.
\end{equation*}
Then it is straightforward to establish the following result.

\vspace*{0.25cm}
\begin{proposition}\label{eq_separable}
Assume that the fertility function $\beta$ satisfies \eqref{specbeta} and that the following conditions hold true
\begin{align}
& \beta(s,y,0)>\mu(s,0),\quad s,y\in[0,m],\, P\in (0,\infty); \quad \int_0^m\exp\left\{-\int_0^s\frac{\mu(y,0)}{\gamma(y,0)}\,\ud y\right\}\,\ud s < m-1,\label{stat1} \\
& \int_0^m\beta_1(s,P)\,\ud s\to 0\quad \text{as}\quad P\to\infty,\quad \text{and}\quad 0<\gamma^*\le\gamma(s,P),\quad s\in[0,m],\, P\in (0,\infty). \label{stat2}
\end{align}
Then model \eqref{equation}-\eqref{initial} admits at least one positive equilibrium solution.
\end{proposition}
\noindent {\bf Proof.} Condition \eqref{stat1} implies
\begin{align}
R(0) & = \int_0^m \exp\left\{-\int_0^s\frac{\mu(y,0)}{\gamma(y,0)}\,\ud y\right\}\int_0^s\frac{\beta_2(s)\beta_1(y,0)}{\gamma(y,0)}\exp\left\{\int_0^y\frac{\mu(z,0)}{\gamma(z,0)}\,\ud z\right\}\,\ud y\,\ud s\nonumber \\
& >  \int_0^m \exp\left\{-\int_0^s\frac{\mu(y,0)}{\gamma(y,0)}\,\ud y\right\}\int_0^s\left(\exp\left\{\int_0^y\frac{\mu(z,0)}{\gamma(z,0)}\,\ud z\right\} \right)'\,\ud y\,\ud s\nonumber \\
& >1.
\end{align}
Condition \eqref{stat2} and the growth behaviour of the functions in \eqref{netrep} imply that 
\begin{equation*}
\displaystyle\lim_{P\to +\infty}R(P)=0,
\end{equation*} 
hence the claim holds true on the grounds of the Intermediate Value Theorem.
\hfill $\Box$

\subsection{The general case}

For a fixed $P\in (0,\infty)$ we define the operator $\mathcal{B}_P$ by
\begin{align}
\mathcal{B}_P\,u= & -\frac{\partial}{\partial s}\left(\gamma(\cdot,P)u\right)-\mu(\cdot,P)u+\int_0^m\beta(\cdot,y,P)u(y)\,\ud y, \nonumber \\
\text{Dom}(\mathcal{B}_P)= & \left\{u\in W^{1,1}(0,m)\,|\,u(0)=0\right\}.\label{operator}
\end{align}
Our goal is to show that there exists a $P_*$ such that the operator $\mathcal{B}_{P_*}$ has eigenvalue $0$ with a corresponding unique positive 
eigenvector. To this end, first we establish that $\mathcal{B}_{P}$ is the generator of a positive semigroup. 
Then we determine conditions that guarantee that it generates an irreducible semigroup. 
We also establish that the governing linear semigroup is eventually compact, which implies that the Spectral Mapping Theorem holds true for the semigroup 
and its generator, and the spectrum of the generator may contain only isolated eigenvalues of finite algebraic multiplicity (see e.g. \cite{NAG}). 
It then follows that the spectral bound is a dominant (real) eigenvalue $\lambda_P$ of geometric multiplicity one with a corresponding  
positive eigenvector \cite[Chapter 9]{CH}. Finally we need to establish conditions which imply that there exist a 
$P^+\in (0,\infty)$ such that the spectral bound $s(\mathcal{B}_{P^+})$ is negative and therefore the dominant eigenvalue 
$\lambda_{P^+}=s(\mathcal{B}_{P^+})$ is also 
negative; and a $P^-\in (0,\infty)$ such that this dominant eigenvalue $\lambda_{P^-}=s(\mathcal{B}_{P^-})$ is positive. 
Then it follows from standard perturbation results on eigenvalues (see e.g. \cite{K}) that there exists a zero eigenvalue.  
A similar strategy was employed in \cite{CS} to establish the existence and uniqueness of an equilibrium solution of a cyclin structured cell population model.  

\vspace*{0.25cm}
\begin{lemma}
For every $P\in (0,\infty)$ the semigroup $\mathcal{T}(t)$ generated by the operator $\mathcal{B}_P$ is positive.
\end{lemma}
\noindent {\bf Proof.} We rewrite \eqref{operator} as, $\mathcal{B}_P=\mathcal{A}_P+\mathcal{C}_P$, where 
\begin{align}
\mathcal{A}_P\,u= & -\frac{\partial}{\partial s}\left(\gamma(\cdot,P)u\right)-\mu(\cdot,P)u \nonumber \\
\text{Dom}(\mathcal{A}_P)= & \left\{u\in W^{1,1}(0,m)\,|\,u(0)=0\right\}, \nonumber \\
\mathcal{C}_P\,u= & \int_0^m\beta(\cdot,y,P)u(y)\,\ud y, \nonumber \\
\text{Dom}(\mathcal{C}_P)= & \,L^1(0,m).\label{operator2}
\end{align}
For $0\le f\in L^1(0,m)$ the solution of the resolvent equation
\begin{equation*}
(\lambda\mathcal{I}-\mathcal{A}_P)u=f,
\end{equation*}
is
\begin{equation*}
u(s)=\int_0^s\exp\left\{-\int_y^s\frac{\lambda+\gamma_s(\sigma,P_*)+\mu(\sigma,P_*)}{\gamma(\sigma,P_*)}\,\ud \sigma\right\}\frac{f(y)}{\gamma(y,P_*)}\,\ud y.
\end{equation*}
This shows that the resolvent operator $\mathcal{R}(\lambda,\mathcal{A}_P)$ is a positive bounded operator, hence $\mathcal{A}_P$ ge\-ne\-rates a 
positive semigroup. Since $\mathcal{C}_P$ is a positive and bounded operator, the statement follows.
\hfill $\Box$
\vspace*{0.25cm}
\begin{lemma}\label{event_compact}
The linear semigroup $\mathcal{T}(t)$ generated by the operator $\mathcal{B}_{P}$ is eventually compact.
\end{lemma}
\noindent {\bf Proof.}  We note that $\mathcal{A}_P$ generates a nilpotent semigroup, while it is easily shown that $\mathcal{C}_P$ is a compact operator 
if conditions \eqref{assumptions} hold true. (For more details see also Theorem \ref{reg}.) \hfill $\Box$

\vspace*{0.25cm}
\begin{lemma}\label{irred}
Assume that for every $P\in (0,\infty)$ there exists an $\varepsilon_0>0$ such that for all $0<\varepsilon\le\varepsilon_0$
\begin{equation}\label{assume_birth}
 \int_0^\varepsilon \int_{m-\varepsilon}^m  \beta(s,y,P)\,\ud y\,\ud s > 0.
\end{equation}
Then the linear semigroup $\mathcal{T}(t)$ generated by the operator $\mathcal{B}_{P}$ is irreducible.
\end{lemma}
\noindent {\bf Proof.}
We only need to show that under condition \eqref{assume_birth} for every $p_0\in L^1_+(0,m)$ there exists a $t_0$ such that 
\begin{equation*}
 supp\, \mathcal{T}(t_0)p_0=[0,m],
\end{equation*}
for all $t\ge t_0$. Since $\gamma>0$, there exists $t_*$ such that 
\begin{equation*}
supp \,\mathcal{T}(t) p_0 \cap supp\,\beta(s,\,\cdot\,)\neq \emptyset
\end{equation*}
for every $t_*\le t$ and every $s\in(0,\varepsilon]$. By assumption \eqref{assume_birth}, \mbox{$\mathcal{T}(t) p_0(s)>0$} for $t_*\le t$ and $s\in(0,\varepsilon]$.  After this, eventually the support of the solution $\mathcal{T}(t_0)p_0$ will cover the entire size space $[0,m]$. \hfill $\Box$

\vspace*{0.25cm}
\begin{lemma}\label{equilibrium}
Assume that there exists a $\beta^-(s,y,P)=\beta^-_1(s,P)\beta^-_2(y)$ and a $P^-\in (0,\infty)$ 
such that 
\begin{equation}
\beta^-_1(s,P^-)\beta^-_2(y)\le\beta(s,y,P^-),\quad s,\,y\in [0,m],\label{1pos}
\end{equation} 
and
\begin{equation}
\int_0^m\int_0^s\frac{\beta^-_1(y,P^-)\beta^-_2(s)}{\gamma(y,P^-)}\exp\left\{-\int_y^s\frac{\gamma_s(z,P^-)+\mu(z,P^-)}{\gamma(z,P^-)}\,\ud z\right\}\,\ud y\,\ud s>1,\label{eqcond1}
\end{equation}
and a $\beta^+(s,y,P)=\beta^+_1(s,P)\beta^+_2(y)$ and a $P^+\in (0,\infty)$ such that 
\begin{equation}
\beta(s,y,P^+)\le\beta^+_1(s,P^+)\beta^+_2(y),\label{2pos}
\end{equation}
and
\begin{equation}
\int_0^m\int_0^s\frac{\beta^+_1(y,P^+)\beta^-_2(s)}{\gamma(y,P^+)}\exp\left\{-\int_y^s\frac{\gamma_s(z,P^+)+\mu(z,P^+)}{\gamma(z,P^+)}\,\ud z\right\}\,\ud y\,\ud s<1.\label{eqcond2}
\end{equation}
Then the operator $\mathcal{B}_{P^-}$ has a dominant real eigenvalue $\lambda_{P^-}>0$ and the operator $\mathcal{B}_{P^+}$ 
has a dominant real eigenvalue $\lambda_{P^+}<0$, with corresponding positive eigenvectors.
\end{lemma}
\noindent {\bf Proof.} 
First assume that there exists a $\beta^-(s,y,P)=\beta^-_1(s,P)\beta^-_2(y)$ and a $P^-$ such that conditions \eqref{1pos} and \eqref{eqcond1} hold true. 
Let $\mathcal{B}_{P^-}^-$ denote the operator that corresponds to the fertility $\beta^-$ and the constant $P^-$.  The solution of the eigenvalue problem
\begin{equation}
\mathcal{B}_{P^-}^-u=\lambda u,\quad\quad u(0)=0
\end{equation}
is 
\begin{equation}
u(s)=\int_0^m\beta^-_2(s)u(s)\,\ud s\int_0^s\frac{\beta^-_1(y,P^-)}{\gamma(y,P^-)}\exp\left\{-\int_y^s\frac{\lambda+\gamma_s(z,P^-)+\mu(z,P^-)}{\gamma(z,P^-)}\,\ud z\right\}\,\ud y.\label{eq1}
\end{equation}
We multiply equation \eqref{eq1} by $\beta^-_2$ and integrate from $0$ to $m$ to arrive at the characteristic equation
\begin{equation}
1=\int_0^m \beta^-_2(s)\int_0^s\frac{\beta^-_1(y,P^-)}{\gamma(y,P^-)}\exp\left\{-\int_y^s\frac{\lambda+\gamma_s(z,P^-)+\mu(z,P^-)}{\gamma(z,P^-)}\,\ud z\right\}\,\ud y\,\ud s.\label{eq2}
\end{equation}
Equation \eqref{eq2} admits a unique dominant real solution $\lambda_{P^-}^->0$ if condition \eqref{eqcond1} holds true.
Since $\mathcal{B}_{P^-}^-$ is a generator of a positive semigroup and $(\mathcal{B}_{P^-}-\mathcal{B}_{P^-}^-)$ is a positive (and bounded) operator 
by condition \eqref{1pos}, it follows that $\mathcal{B}_{P^-}$ has a dominant real eigenvalue $\lambda_{P^-}\ge\lambda_{P^-}^->0$, 
see e.g.~\cite[Corollary VI.1.11]{NAG}.

In a similar way, let us assume that there exists a $\beta^+(s,y,P)=\beta^+_1(s,P)\beta^+_2(y)$ and a $P^+$ such that condition \eqref{2pos} and \eqref{eqcond2} hold true. 
Let $\mathcal{B}_{P^+}^+$ denote the operator which corresponds to the fertility $\beta^+$ and the constant $P^+$. 
The solution of the eigenvalue problem
\begin{equation}
\mathcal{B}_{P^+}^+u=\lambda u,\quad\quad u(0)=0
\end{equation}
is now
\begin{equation}
u(s)=\int_0^m\beta^+_2(s)u(s)\,\ud s\int_0^s\frac{\beta^+_1(y,P^+)}{\gamma(y,P^+)}\exp\left\{-\int_y^s\frac{\lambda+\gamma_s(z,P^+)+\mu(z,P^+)}{\gamma(z,P^+)}\,\ud z\right\}\,\ud y.\label{eq3}
\end{equation}
We multiply equation \eqref{eq3} by $\beta^+_2$ and integrate from $0$ to $m$ to arrive at the characteristic equation
\begin{equation}
1=\int_0^m \beta^+_2(s)\int_0^s\frac{\beta^+_1(y,P^+)}{\gamma(y,P^+)}\exp\left\{-\int_y^s\frac{\lambda+\gamma_s(z,P^+)+\mu(z,P^+)}{\gamma(z,P^+)}\,\ud z\right\}\,\ud y\,\ud s.\label{eq4}
\end{equation}
Equation \eqref{eq4} admits a unique dominant real solution $\lambda_{P^+}^+<0$ if condition \eqref{eqcond2} holds true.
Since $\mathcal{B}_{P^+}$ is a generator of a positive semigroup and $(\mathcal{B}_{P^+}^+-\mathcal{B}_{P^+})$ is a positive operator 
by condition \eqref{2pos}, it follows that $\mathcal{B}_{P^+}$ has a dominant real eigenvalue $\lambda_{P^+}\le\lambda_{P^+}^+<0$. 

In both cases, the positivity of the corresponding eigenvector follows from the irreducibility of the semigroup $\mathcal{T}(t)$, 
see \cite[Theorem 9.11]{CH}.
\hfill $\Box$

\vspace*{0.25cm}
\begin{theorem}\label{equilibrium2}
Assume that conditions \eqref{assume_birth}, \eqref{1pos}-\eqref{eqcond2} are satisfied. Then system \eqref{equation}-\eqref{initial} admits at least one 
positive equilibrium solution.
\end{theorem}
\noindent {\bf Proof.}  Let $P^*>0$ be such that $s(\mathcal{B}_P^*)=0$. Then, since the spectrum consists only of isolated eigenvalues we have 
$\lambda_{P^*}=s(\mathcal{B}_{P^*})=0$ and there exists a corresponding positive eigenvector $p_*$. 
Then $\frac{P^*}{||p_*||_1} p_*$  is the desired equilibrium solution with total population size $P^*$.  \hfill $\Box$

\vspace*{0.5cm}
\setcounter{equation}{0}
\section{The linearised semigroup and its regularity}\label{semigroup}

Here, when we use the term `linearised semigroup', we refer to the linear semigroup governing the linearised system. 
However, since it was proved in \cite{CS2} that model \eqref{equation}-\eqref{initial} is well-posed, there exists a semigroup of nonlinear operators  
$\Sigma(t)_{t\ge 0}$ defined via $\Sigma(t)p(s,0)=p(s,t)$. It was proven in \cite{D2} that if the nonlinearities are 
smooth enough (namely, the vital rates are differentiable) then this nonlinear semigroup $\Sigma(t)$ is Frech\'{e}t differentiable and the Frech\'{e}t 
derivative around an equilibrium solution $p_*$ defines a semigroup of bounded linear operators. In this section we will establish the existence of this 
semigroup and at the same time arrive at a condition which guarantees that it is positive.

Given a positive stationary solution $p_*$ of system \eqref{equation}-\eqref{initial}, we introduce the perturbation $u=u(s,t)$ of $p$ by 
making the ansatz $p=u+p_*$. A Taylor series expansion of the vital rates gives the linearised problem (see e.g. \cite{FH})
\begin{align}
u_t(s,t) &  =  -\gamma(s,P_*)\,u_s(s,t)-\left(\gamma_s(s,P_*)+ \mu(s,P_*)\right)\,u(s,t)\nonumber\\
&\quad -\left(\gamma_{sP}(s,P_*)\,p_*(s)+\mu_P(s,P_*)\,p_*(s)+\gamma_P(s,P_*)\,p_*�^\prime(s)\right)\,U(t) \nonumber \\
&\quad +\int_0^m u(y,t)\left(\beta(s,y,P_*)+\int_0^m\beta_P(s,z,P_*)p_*(z)\,\ud z\right)\,\ud y,\label{lin1}\\
\gamma(0,P_*) u(0,t) & =0\label{lin2}
\end{align}
where we have set
\begin{equation}
    U(t)=\int_0^m u(s,t)\,\ud s.
\end{equation}
Eqs.~\eqref{lin1}--\eqref{lin2} are accompanied by the initial condition
\begin{equation}\label{lin3}
    u(s,0)=u_0(s).
\end{equation}
Our first objective is to establish conditions which guarantee that the linearised system is governed by a positive semigroup.  
To this end, we cast the linearised system \eqref{lin1}-\eqref{lin3} in the form of an abstract Cauchy problem on the state space 
$\mathcal{X}=L^1(0,m)$ as follows
\begin{equation}\label{abstr}
    \frac{d}{dt}\, u = \left({\mathcal A} + {\mathcal B}+ {\mathcal C}+{\mathcal D}\right)\,u,\quad    u(0)=u_0,
\end{equation}
where
\begin{align}
& {\mathcal A} u = -\gamma(\cdot,P_*)\,u_s\quad \text{with domain}\quad\text{Dom}({\mathcal A})=\left\{u\in W^{1,1}(0,m)\,|\,u(0)=0\right\},\\
& {\mathcal B} u = -\left(\gamma_s(\cdot,P_*)+ \mu(\cdot,P_*)\right)\,u\quad \text{on ${\mathcal X}$,}\\
& {\mathcal C} u = -\left(\gamma_{sP}(\cdot,P_*)\,p_*+\mu_P(\cdot,P_*)\,p_*+\gamma_P(\cdot,P_*)\,p_*^\prime\right)\,\int_0^m u(s)\,\ud s \nonumber \\
& \quad\, =-\rho_*(\cdot)\int_0^mu(s)\,\ud s\quad \text{on ${\mathcal X}$,} \label{rho} \\
& {\mathcal D} u=\int_0^m u(y)\left(\beta(\cdot,y,P_*)+\int_0^m\beta_P(\cdot,z,P_*)p_*(z)\,\ud z\right)\,\ud y\quad \text{on ${\mathcal X}$},
\end{align}
where $\rho_*$ is defined via equation \eqref{rho}. Our aim is to establish that the linear operator $\mathcal{A+B+C+D}$ is a generator of a quasicontraction semigroup. To this end first we recall (see e.g. \cite{AGG,CH,NAG}) some basic concepts from the theory of linear operators acting on Banach spaces. Let $\mathcal{O}$ be a linear operator defined on the real Banach space $\mathcal{Y}$ with norm $||.||$. $\mathcal{O}$ is called dissipative if for every $\lambda>0$ and $x\in \text{Dom}(\mathcal{O})$,
\begin{equation*}
||(\mathcal{I}-\lambda\mathcal{O})x||\ge ||x||.
\end{equation*}
Furthermore, a function $f\,:\,\mathcal{Y}\to\mathbf{R}$ is called sublinear if
\begin{align*}
& f(x+y)\le f(x)+f(y),\quad x,y\in\mathcal{Y}\\
& f(\lambda x)=\lambda f(x),\quad \lambda\ge 0,\quad x\in\mathcal{Y}.
\end{align*}
If also $f(x)+f(-x)>0$ holds true for $x\ne 0$ then $f$ is called a half-norm on $\mathcal{Y}$. The linear operator $\mathcal{O}$ is called $f$-dissipative if
\begin{equation*}
f(x)\le f(x-\lambda \mathcal{O}x),\quad \lambda\ge 0,\quad x\in \text{Dom}(\mathcal{O}).
\end{equation*}
An operator $\mathcal{O}$ which is $p$-dissipative with respect to the half norm
\begin{equation*}
p(x)=||x^+||,
\end{equation*}
is called dispersive, where $x^+=x\vee 0$ (and $x^-=(-x)^+$). Finally a $C_0$ semigroup $\left\{\mathcal{T}(t)\right\}_{t\ge 0}$ is called quasicontractive if
\begin{equation*}
||\mathcal{T}(t)||\le e^{\omega t},\quad t\ge 0,
\end{equation*}
for some $\omega\in\mathbf{R}$, and it is called contractive if $\omega\le 0$.
We recall the following characterization theorem from \cite{CH}.
\begin{theorem}\label{dispersive}
Let $\mathcal{Y}$ be a Banach lattice and let $\mathcal{O}\,:\,\text{Dom}(\mathcal{O})\to\mathcal{Y}$ be a linear operator. Then, the following statements are equivalent.
\begin{enumerate}
\item[(i)] $\mathcal{O}$ is the generator of a positive contraction semigroup.
\item[(ii)] $\mathcal{O}$ is densely defined, $\text{Rg}(\lambda\mathcal{I}-\mathcal{O})=\mathcal{Y}$ for some $\lambda>0$, and $\mathcal{O}$ is dispersive.
\end{enumerate}
\end{theorem}
We also recall that $\mathcal{O}$ is dispersive if for every $x\in\text{Dom}(\mathcal{O})$ there exists $\phi\in\mathcal{Y}^*$ with $0\le\phi$,  $||\phi||\le 1$ and $(x,\phi)=||x^+||$ such that $(\mathcal{O}x,\phi)\le 0$, where $(\cdot\,,\cdot)$ is the natural pairing between elements of $\mathcal{Y}$ and its dual $\mathcal{Y}^*$.

\vspace*{0.25cm}
\begin{theorem}\label{posgen}
The operator ${\mathcal A} + {\mathcal B}+ {\mathcal C}+{\mathcal D}$ generates a positive strongly continuous ($C_0$ for short) quasicontraction semigroup $\{{\mathcal T}(t)\}_{t\geq 0}$ of bounded linear operators on ${\mathcal X}$ if the following condition holds true
\begin{equation}\label{poscond1}
\rho_*(s)\le\beta(s,y,P_*)+\int_0^m\beta_P(s,y,P_*)p_*(y)\,\ud y,\quad s,y\in [0,m],
\end{equation}
where $\rho_*$ is defined via equation \eqref{rho}.
\end{theorem}

\noindent {\bf Proof.} Our aim is to apply the previous characterization theorem for the perturbed operator $\mathcal{A+B+C+D}-\omega\mathcal{I}$, for some $\omega\in\mathbf{R}$. To this end, for every $u\in\text{Dom}(\mathcal{A+B+C+D}-\omega\mathcal{I})$ we define $\phi_u\in\mathcal{X}^*$ by
\begin{align}
& \phi_u(s)=\frac{u^+(s)}{|u(s)|}, \quad s\in[0,m],\quad u(s)\ne 0,
\end{align}
if $u(s)=0$ then let $\phi_u(s)=0$. Then
\begin{equation*}
||\phi_u||_\infty\le 1,
\end{equation*}
and clearly
\begin{equation*}
(u,\phi_u)=\int_0^mu(s)\phi_u(s)\,\ud s=||u^+||_1.
\end{equation*}
Making use of condition \eqref{poscond1} we obtain the following estimate.
\begin{align}
& ((\mathcal{A+B+C+D} -\omega\mathcal{I})u,\phi_u) \nonumber \\
& = -\int_0^m \mathbf{1}_{u^+}(s)\,\big(\gamma(s,P_*)u(s)\big)_s\,\ud s-\int_0^m\mathbf{1}_{u^+}(s)\,\mu(s,P_*)u(s)\,\ud s-\int_0^m\mathbf{1}_{u^+}(s)\,\omega\,u(s)\,\ud s \nonumber \\
& \quad\, +\int_0^m \mathbf{1}_{u^+}(s)\,\int_0^m u(y)\left(\beta(s,y,P_*)+\int_0^m\beta_P(s,z,P_*)p_*(z)\,\ud z-\rho_*(s)\right)\,\ud y\,\ud s \nonumber \\
& \le -\int_0^m \mathbf{1}_{u^+}(s)\,\big(\gamma(s,P_*)u(s)\big)_s\,\ud s-\omega||u^+||_1-\inf_{s\in [0,m]}\mu(s,P_*)\,||u^+||_1 \nonumber \\
& \quad\, +||u^+||_1\left|\left|\sup_{y\in[0,m]}\left(\beta(s,y,P_*)+\int_0^m\beta_P(s,z,P_*)p_*(z)\,\ud z-\rho_*(s)\right)\right|\right|_\infty\nonumber \\
& \le -\omega||u^+||_1-(\gamma(m,P_*)u(m))\mathbf{1}_{u^+}(m)\nonumber \\
& \quad\, +||u^+||_1\left|\left|\sup_{y\in[0,m]}\left(\beta(s,y,P_*)+\int_0^m\beta_P(s,z,P_*)p_*(z)\,\ud z-\rho_*(s)\right)\right|\right|_\infty\nonumber \\
& \le 0, \label{estimates}
\end{align}
for some $\omega\in\mathbf{R}$ large enough, hence the operator $\mathcal{A+B+C+D}-\omega\mathcal{I}$ is dispersive.
The operator $\mathcal{A+B+C+D}-\omega\mathcal{I}$ is clearly densely defined.  We observe that the equation
\begin{equation}
    (\lambda I-\mathcal{A})\,u=h
\end{equation}
for $h\in\mathcal{X}$ and $\lambda>0$ sufficiently large has a unique solution $u\in \text{Dom}({\mathcal A})$, given by
\begin{align}
u(s)=\exp\left\{-\int_0^s\frac{\lambda}{\gamma(y,P_*)}\,\ud y\right\}\int_0^s\exp\left\{\int_0^y\frac{\lambda}{\gamma(z,P_*)}\,\ud z\right\}\frac{h(y)}{\gamma(y,P_*)}\,\ud y.\label{udef1}
\end{align}
The fact that $u\in \text{Dom}({\mathcal A})$ is well defined by \eqref{udef1} follows from
\begin{align*}
|u'(s)| & \le \left|\frac{h(s)}{\gamma(s,P_*)}\right|+\frac{\lambda}{\gamma(s,P_*)}\int_0^m\exp\left\{-\int_y^s\frac{\lambda}{\gamma(z,P_*)}\,\ud z\right\}\frac{|h(y)|}{\gamma(y,P_*)}\,\ud y \nonumber \\
& \le \left|\frac{h(s)}{\gamma(s,P_*)}\right|+M_\lambda,
\end{align*}
for $\lambda$ large enough for some $M_\lambda<\infty$, that is $u\in W^{1,1}(0,m)$. Since $\mathcal{B+C+D}-\omega\mathcal{I}$ is bounded, the range condition is satisfied. Theorem \ref{dispersive} gives that $\mathcal{A+B+C+D}-\omega\mathcal{I}$ is a generator of a positive contraction semigroup. Since the operator $\omega\mathcal{I}$ is positive (clearly if the dispersivity estimate holds true with an $\omega<0$ then it holds true with any other $\omega^*>\omega$) a well-known perturbation result (see e.g. \cite{NAG}) yields that $\mathcal{A+B+C+D}$ is a generator of a positive quasicontraction semigroup $\mathcal{T}$ which obeys
\begin{equation*}
\|\mathcal{T}(t)\|\le e^{\omega t},\quad t\geq 0.
\end{equation*}
\hfill $\Box$

\vspace*{0.25cm}
\begin{remark}\label{remark1}
The proof of Theorem \ref{dispersive} shows that if
\begin{align*}
& \inf_{s\in[0,m]}\mu(s,P_*)>\left|\left|\sup_{y\in[0,m]}\left(\beta(s,y,P_*)+\int_0^m\beta_P(s,z,P_*)p_*(z)\,\ud z-\rho_*(s)\right)\right|\right|_\infty
\end{align*}
holds, then the growth bound $\omega_0$ of the semigroup is negative, hence the semigroup $\{\mathcal{T}(t)\}_{t\ge 0}$ is uniformly exponentially stable (see e.g. \cite{NAG}), i.e. the equilibrium $p_*$ is locally asymptotically stable.
\end{remark}

\vspace*{0.25cm}
\begin{remark}
We note that the operator $\mathcal{A+B+C+D}$ is in general a generator of a $C_0$ quasicontraction (but not positive) semigroup. The proof of this would utilise the Lumer-Phillips Theorem (see e.g. \cite{AGG,CH,NAG}) and goes along similar lines, obtaining a dissipativity estimate in terms of $u$ rather than $u^+$, see e.g. \cite{FH2}. This implies that the linearised problem \eqref{lin1}-\eqref{lin2} is well-posed.
\end{remark}

\vspace*{0.25cm}
\begin{remark}
Note that if $\beta=\beta(s,y),\,\mu=\mu(s),\,\gamma=\gamma(s)$, i.e. model \eqref{equation}-\eqref{initial} is a linear one, then the biologically relevant conditions $\mu,\beta\ge 0$ and $\gamma>0$ imply that it is governed by a positive quasicontraction semigroup.
\end{remark}

\vspace*{0.25cm}
\begin{theorem}\label{reg}
The semigroup $\{{\mathcal T}(t)\}_{t\geq 0}$ generated by the operator ${\mathcal A} + {\mathcal B}+ {\mathcal C}+ {\mathcal D}$ is eventually compact.
\end{theorem}
\noindent {\bf Proof.}
${\mathcal C}$ is a rank-one operator. Hence it is compact on ${\mathcal X}=L^1(0,m)$. $\mathcal{D}$ is linear and bounded. Hence in view of the Fr\'{e}chet-Kolmogorov compactness criterion in $L^p$ we need to show that
\begin{equation*}
\lim_{t\to 0}\int_0^m\left|\mathcal{D}u(t+s)-\mathcal{D}u(s)\right|\,\ud s=0,\quad \text{uniformly in}\,\,u,
\end{equation*}
for $u\in B$, where $B$ is the unit sphere of $L^1(0,m)$. But this follows from the regularity assumptions we made on $\beta$ based on the following estimate
\begin{align*}
& |\mathcal{D}u(s_1)-\mathcal{D}u(s_2)|\le ||u||_1 \\
& \times\left|\left|\beta(s_1,y,P_*)+\int_0^m\beta_P(s_1,z,P_*)p_*(z)\,\ud z-\beta(s_2,y,P_*)-\int_0^m\beta_P(s_2,z,P_*)p_*(z)\,\ud z\right|\right|_\infty.
\end{align*}
Therefore, it suffices to investigate the operator ${\mathcal A}+{\mathcal B}$. To this end, we note that the abstract differential equation
\begin{equation}
\frac{d}{dt}\, u = ({\mathcal A}+{\mathcal B})\,u
\end{equation}
corresponds to the partial differential equation
\begin{equation}\label{reduc}
u_t(s,t)+ \gamma(s,P_*)\,u_s(s,t)+\left(\gamma_s(s,P_*)+ \mu(s,P_*)\right)\,u(s,t) = 0,
\end{equation}
subject to the boundary condition \eqref{lin2}. We solve easily equation \eqref{reduc} using the method of characteristics. For $t>\Gamma(m)$ we arrive at
\begin{equation}
u(s,t)=u(0,t-\Gamma(s))\exp\left\{-\int_0^s\frac{\gamma_s(y,P_*)+\mu(y,P_*)}{\gamma(y,P_*)}\,\ud y\right\}=0,
\end{equation}
where
\begin{equation*}
\Gamma(s)=\int_0^s\frac{1}{\gamma(y,P_*)}\,\ud y.
\end{equation*}
This means that the semigroup $\mathcal{T}(t)$ generated by $\mathcal{A+B}$ is nilpotent. In particular it is compact for $t>\Gamma(m)$ and the claim follows.
\hfill $\Box$

\vspace*{0.25cm}
\begin{remark}
Theorem \ref{reg} implies that the Spectral Mapping Theorem holds true for the semigroup $\left\{\mathcal{T}(t)\right\}_{t\ge 0}$ with generator $\mathcal{A+B+C+D}$ and that the spectrum $\sigma(\mathcal{A+B+C+D})$ contains only isolated eigenvalues of finite multiplicity (see e.g. \cite{NAG}).
\end{remark}

\section{(In)\,Stability}

Here, we consider the stability of positive equilibrium solutions by studying the point spectrum of the linearised operator $\mathcal{A+B+C+D}$. The main difficulty is  that the eigenvalue equation
\begin{equation*}
(\mathcal{A+B+C+D}-\mathcal{I})\lambda=0,
\end{equation*}
cannot be solved explicitly, since in general, the operator $\mathcal{D}$ has infinite rank. We encountered this problem previously with hierarchical size-structured population models \cite{FH2,FH3}. In \cite{FH2} and \cite{FH3} we used the dissipativity approach, 
presented in the previous section, to establish conditions which guarantee that the spectral bound of the linearised semigroup is negative. 
However, as we can see from Remark \ref{remark1} this approach gives a rather restrictive stability condition. Therefore, here we devise a different 
approach, which uses positive perturbation arguments.

\vspace*{0.25cm}
\begin{theorem}\label{instab}
Assume that there exists an $\varepsilon>0$ such that
\begin{equation}\label{instabcond1}
\beta(s,y,P_*)-\rho_*(s)-\varepsilon+\int_0^m\beta_P(s,y,P_*)p_*(y)\,\,d y\ge 0,\quad s,y\in [0,m],
\end{equation}
and 
\begin{equation}\label{instabcond2}
\varepsilon\int_0^m\exp\left\{-\int_0^s\frac{\gamma_s(\sigma, P_*)+\mu(\sigma,P_*)}{\gamma(\sigma,P_*)}\,\ud\sigma  \right\}\int_0^s\frac{\exp\left\{\int_0^y\frac{\gamma_s(\sigma,P_*)+\mu(\sigma,P_*)}{\gamma(\sigma,P_*)}\,\ud\sigma \right\}}{\gamma(y,P_*)}\,\ud y\,\ud s>1.
\end{equation}
Then the stationary solution $p_*(s)$ of model \eqref{equation}-\eqref{initial} is linearly unstable.
\end{theorem}
\noindent {\bf Proof.} Let $\varepsilon>0$, and define the operator $\mathcal{F}_\varepsilon$ on $\mathcal{X}$ as
\begin{equation*}
\Fe u = \varepsilon \int_0^m u(s)\,\ud s =  \varepsilon \bar{u}.
\end{equation*}
We first find the solution of the eigenvalue equation
\begin{equation*}
(\A+\B+\Fe) u = \lambda u
\end{equation*}
as
\begin{align}
u(s) = &  \,\varepsilon\, \bar{u}\, \exp\left\{-\int_0^s\frac{\lambda+\gamma_s(\sigma,P_*)+\mu(\sigma,P_*)}{\gamma(\sigma,P_*)}\,\ud\sigma  \right\} \nonumber \\ 
& \times\int_0^s\frac{1}{\gamma(y,P_*)} \exp\left\{\int_0^y\frac{\lambda+\gamma_s(\sigma,P_*)+\mu(\sigma,P_*)}{\gamma(\sigma,P_*)}\,\ud\sigma  \right\}\,\ud y.\label{eigv}
\end{align}
Next we integrate the solution \eqref{eigv} over $[0,m]$ to obtain
\begin{align}
\bar{u} = & \, \varepsilon\, \bar{u}\int_0^m\left[ \exp\left\{-\int_0^s\frac{\lambda+\gamma_s(\sigma,P_*)+\mu(\sigma,P_*)}{\gamma(\sigma,P_*)}\,\ud\sigma  \right\} \right. \nonumber \\
 & \left. \quad\quad\times \int_0^s\frac{1}{\gamma(y,P_*)} \exp\left\{\int_0^y\frac{\lambda+\gamma_s(\sigma,P_*)+\mu(\sigma,P_*)}{\gamma(\sigma,P_*)}\,\ud\sigma  \right\}\,\ud y\right]\,\ud s.\label{eigv2}
\end{align}
We note that, if $\bar{u}=0$ then equation \eqref{eigv} shows that $u(s)\equiv 0$, hence we have a non-trivial eigenvector if and only if 
$\bar{u}\ne 0$ and $\lambda$ satisfies the following characteristic equation
\begin{align}
1 =K(\lambda)\eqdef  & \varepsilon \int_0^m\left[ \exp\left\{-\int_0^s\frac{\lambda+\gamma_s(\sigma,P_*)+\mu(\sigma,P_*)}{\gamma(\sigma,P_*)}\,\ud\sigma  \right\} \right. \nonumber \\
 & \left. \quad\quad \times\int_0^s\frac{1}{\gamma(y,P_*)} \exp\left\{\int_0^y\frac{\lambda+\gamma_s(\sigma,P_*)+\mu(\sigma,P_*)}{\gamma(\sigma,P_*)}\,\ud\sigma \right\}\,\ud y\right]\,\ud s. \label{char1}
\end{align}
It is easily shown that
\begin{equation*}
\lim_{\lambda\ra +\infty} K(\lambda)  = 0,
\end{equation*}
therefore it follows from condition \eqref{instabcond2}, on the grounds of the Intermediate Value Theorem, that equation \eqref{char1} has a positive (real) solution. Hence we have
\begin{equation*}
0<s(\mathcal{A+B+F_{\varepsilon}}).
\end{equation*}
Next, for a fixed $0\le f\in\mathcal{X}$, we obtain the solution of the resolvent equation
\begin{equation*}
\left(\lambda\mathcal{I}-(\mathcal{A+B+F_{\varepsilon}})\right)u=f,
\end{equation*}
as
\begin{align}
u(s) = &  \exp\left\{-\int_0^s\frac{\lambda+\gamma_s(\sigma,P_*)+\mu(\sigma,P_*)}{\gamma(\sigma,P_*)}\,\ud\sigma  \right\} \nonumber \\ 
& \times\int_0^s \exp\left\{\int_0^y\frac{\lambda+\gamma_s(\sigma,P_*)+\mu(\sigma,P_*)}{\gamma(\sigma,P_*)}\,\ud\sigma  \right\}\frac{\varepsilon\bar{u}+f(y)}{\gamma(y,P_*)}\,\ud y.\label{res}
\end{align}
We integrate equation \eqref{res} from $0$ to $m$ to obtain
\begin{equation}
\bar{u}=\frac{\int_0^m\exp\left\{-\int_0^s\frac{\lambda+\gamma_s(\sigma,P_*)+\mu(\sigma,P_*)}{\gamma(\sigma,P_*)}\,\ud\sigma  \right\}\int_0^s \exp\left\{\int_0^y\frac{\lambda+\gamma_s(\sigma,P_*)+\mu(\sigma,P_*)}{\gamma(\sigma,P_*)}\,\ud\sigma  \right\}\frac{f(y)}{\gamma(y,P_*)}\,\ud y}{1-\varepsilon\int_0^m\exp\left\{-\int_0^s\frac{\lambda+\gamma_s(\sigma,P_*)+\mu(\sigma,P_*)}{\gamma(\sigma,P_*)}\,\ud\sigma  \right\}\int_0^s \frac{\exp\left\{\int_0^y\frac{\lambda+\gamma_s(\sigma,P_*)+\mu(\sigma,P_*)}{\gamma(\sigma,P_*)}\,\ud\sigma  \right\}}{\gamma(y,P_*)}\,\ud y}
\end{equation}
It follows from the growth behaviour of the exponential function and from assumptions \eqref{assumptions}, that $\bar{u}$ is well-defined and non-negative 
for any $0\le f\in\mathcal{X}$ and $\lambda$ large enough. Hence the resolvent operator 
\begin{equation*}
\mathcal{R}(\lambda,\mathcal{A+B+F_{\varepsilon}})=(\lambda-(\mathcal{A+B+F_{\varepsilon}}))^{-1}
\end{equation*}
is positive, for $\lambda$ large enough, which implies that $\mathcal{A+B+F_{\varepsilon}}$ generates a positive semigroup (see e.g. \cite{NAG}).

Finally, we note that condition \eqref{instabcond1} guarantees that the operator $\mathcal{C+D-F_{\varepsilon}}$ is positive, hence we have for the spectral bound (see e.g. Corollary VI.1.11 in \cite{NAG})
\begin{equation*}
0<s(\A+\B+\Fe)\le s(\A+\B+\Fe+\C+\D-\Fe) =s(\A+\B+\C+\D),
\end{equation*}
and the result follows.
\hfill $\Box$

Next we show that for a separable fertility function we can indeed explicitly characterize the point spectrum of the linearised operator. 
\begin{theorem}\label{eigenvalues}
Assume that $\beta(s,y,P)=\beta_1(s,P)\beta_2(y),\quad s,y\in [0,m],\, P\in(0,\infty)$.
Then for any $\lambda\in \mathbb{C}$, we have $\lambda\in\sigma(\mathcal{A+B+C+D})$ if and only if $\lambda$ satisfies the equation
\begin{equation}\label{chareq}
K_{\beta}\,(\lambda)=\det
\left( \begin{array}{lll}
1+a_1(\lambda) & a_2(\lambda) \\
a_3(\lambda) & 1+a_4(\lambda) \\
\end{array}\right)=0,
\end{equation}
where
\begin{align}
& a_1(\lambda)=-\int_0^mF(\lambda,s,P_*)\int_0^s\frac{g(y)}{F(\lambda,y,P_*)}\,\ud y\,\ud s,\nonumber \\
& a_2(\lambda)=-\int_0^mF(\lambda,s,P_*)\int_0^s\frac{\beta_1(y,P_*)}{\gamma(y,P_*)F(\lambda,y,P_*)}\,\ud y\,\ud s,\nonumber \\
& a_3(\lambda)=-\int_0^m\beta_2(s)F(\lambda,s,P_*)\int_0^s\frac{g(y)}{F(\lambda,y,P_*)}\,\ud y\,\ud s,\nonumber \\
& a_4(\lambda)=-\int_0^m\beta_2(s)F(\lambda,s,P_*)\int_0^s\frac{\beta_1(y,P_*)}{\gamma(y,P_*)F(\lambda,y,P_*)}\,\ud y\,\ud s,\label{chareq2}
\end{align}
and
\begin{align*}
g(s) & =\frac{\beta_{1_P}(s,P_*)\displaystyle\int_0^m\beta_2(y)p_*(y)\,\ud y-\rho_*(s)}{\gamma(s,P_*)},\quad s\in[0,m], \\
F(\lambda,s,P_*) & =\exp\left\{-\int_0^s\frac{\lambda+\gamma_s(y,P_*)+\mu(y,P_*)}{\gamma(y,P_*)}\,\ud y\right\},\quad s\in [0,m].
\end{align*}
\end{theorem}
\noindent {\bf Proof.}
To characterize the point spectrum of $\mathcal{A+B+C+D}$ we  consider the eigenvalue problem
\begin{equation}\label{eigvalue}
(\mathcal{A+B+C+D-\lambda\mathcal{I}})U=0,\quad U(0)=0.
\end{equation}
The solution of \eqref{eigvalue} is found to be
\begin{align}
U(s)= & \overline{U}F(\lambda,s,P_*)\int_0^s\frac{g(y)}{F(\lambda,y,P_*)}\,\ud y+\widetilde{U}F(\lambda,s,P_*)\int_0^s\frac{\beta_1(y,P_*)}{\gamma(y,P_*)F(\lambda,y,P_*)}\,\ud y,\label{eigsol}
\end{align}
where
\begin{equation*}
\overline{U}=\int_0^m U(s)\,\ud s,\quad \widetilde{U}=\int_0^m \beta_2(s)U(s)\,\ud s.
\end{equation*}
We integrate equation \eqref{eigsol} from zero to $m$ and mulitply equation \eqref{eigsol} by $\beta_2(s)$ and then integrate from zero to $m$ to obtain
\begin{align}
& \overline{U}(1+a_1(\lambda))+\widetilde{U}a_2(\lambda)=0,\label{homog1} \\
& \overline{U}a_3(\lambda)+\widetilde{U}(1+a_4(\lambda))=0.\label{homog2}
\end{align}
If $\lambda\in\sigma(\mathcal{A+B+C+D})$ then the eigenvalue equation \eqref{eigvalue} admits a non-trivial solution $U$ hence there exists a non-zero vector $(\overline{U},\widetilde{U})$ which solves equations \eqref{homog1}-\eqref{homog2}. However, if $(\overline{U},\widetilde{U})$ is a non-zero solution of equations \eqref{homog1}-\eqref{homog2} for some $\lambda\in\mathbb{C}$ then \eqref{eigsol} yields a non-trivial solution $U$. This is because the only scenario for $U$ to vanish would yield
\begin{equation*}
\overline{U}F(\lambda,s)\int_0^s\frac{g(y)}{F(\lambda,y)}\,\ud y=-\widetilde{U}F(\lambda,s)\int_0^s\frac{\beta_1(y,P_*)}{\gamma(y,P_*)F(\lambda,y)}\,\ud y,\quad s\in[0,m].
\end{equation*}
This however, together with equations \eqref{homog1}-\eqref{homog2} would imply $\overline{U}=\widetilde{U}=0$, a contradiction, hence the proof is completed.
\hfill $\Box$

\vspace*{0.25cm}
\begin{theorem}
Assume that condition \eqref{poscond1} holds true for some stationary solution $p_*$.  Moreover, assume that there exists a function $\widetilde{\beta}(s,y,P)=\beta_1(s,P)\beta_2(y)$ such that $\beta(s,y,P_*)\le\widetilde{\beta}(s,y,P_*)$ for $s,y \in [0,m]$ and the characteristic equation $K_{\widetilde{\beta}}\,(\lambda)=0$ does not have a solution with non-negative real part. Then the equilibrium solution $p_*$ is linearly asymptotically stable. 
\end{theorem}

\noindent {\bf Proof.}  We need to establish that the spectral bound of the linearised operator $\mathcal{A+B+C+D}$ is negative. 
To this end, we rewrite the operator $\mathcal{D}$ as a sum of two operators, namely $\mathcal{D}=\mathcal{G}+\mathcal{H}_\beta$, where
\begin{align*}
\mathcal{G}u & =\int_0^m u(y)\,\ud y\int_0^m\beta_P(\cdot,z,P_*)p_*(z)\,\ud z, \quad \text{on} \quad \mathcal{X},\\
\mathcal{H}_\beta u & =\int_0^mu(y)\beta(\cdot,y,P_*)\,\ud y, \quad \text{on} \quad \mathcal{X}.\\
\end{align*}
Condition \eqref{poscond1} guarantees that $\mathcal{A+B+C+G}+\mathcal{H}_\beta$ is a generator of a positive semigroup, while the eventual compactness of the linearised semigroup assures that the spectrum of $\mathcal{A+B+C+G}+\mathcal{H}_{\widetilde{\beta}}$ contains only eigenvalues and that the Spectral Mapping Theorem holds true. Since $\mathcal{H}_{\widetilde{\beta}}-\mathcal{H}_\beta$ is a positive and bounded operator we have
\begin{equation}\label{spectralineq}
s(\mathcal{A+B+C+G}+\mathcal{H}_\beta)\le s(\mathcal{A+B+C+G}+\mathcal{H}_{\beta}+\mathcal{H}_{\widetilde{\beta}}-\mathcal{H}_{\beta})=s(\mathcal{A+B+C+G}+\mathcal{H}_{\widetilde{\beta}})<0,
\end{equation} 
and the proof is completed.
\hfill $\Box$

\begin{example} As we can see from equations \eqref{chareq}-\eqref{chareq2} the characteristic function $K_{\widetilde{\beta}}(\lambda)$ is rather complicated, in general. Therefore, here we only present a special case when it is straightforward to establish that the point spectrum of the linear operator $\mathcal{A+B+C+G}+\mathcal{H}_{\widetilde{\beta}}$ does not contain any element with non-negative real part. In particular, we make the following specific assumption
\begin{equation*}
\beta_2(\cdot)\equiv \beta_2.
\end{equation*}
In this case we can cast the characteristic equation \eqref{chareq} in the simple form
\begin{equation}\label{reducedchar3}
\int_0^m\int_0^s\exp\left\{-\int_y^s\frac{\lambda+\gamma_s(r,P_*)+\mu(r,P_*)}{
\gamma(r,P_*)}\,
\ud r\right\}\left(\frac{g(y)\gamma(y,P_*)+\beta_1(y,P_*)\beta_2}{\gamma(y,P_*)}
\right)\,\ud y\,\ud s=1.
\end{equation}
We note that, if 
\begin{equation*}
g(y)\gamma(y,P_*)+\beta_1(y,P_*)\beta_2\ge 0,\quad\quad y\in [0,m],
\end{equation*}
which is equivalent to the positivity condition \eqref{poscond1}, then equation \eqref{reducedchar3} admits a dominant unique (real) solution. 
On the other hand, it is easily shown that this dominant eigenvalue is negative if 
\begin{equation}
\int_0^m\int_0^s\exp\left\{-\int_y^s\frac{\gamma_s(r,P_*)+\mu(r,P_*)}{\gamma(r,
P_*)}\,
\ud r\right\}\left(\frac{g(y)\gamma(y,P_*)+\beta_1(y,P_*)\beta_2}{\gamma(y,P_*)}
\right)\,\ud y\,\ud s<1.\label{stabcondition}
\end{equation} 
It is easy to see, making use of equation \eqref{netrep}, that \eqref{stabcondition} is satisfied if
\begin{equation*}
\int_0^m\frac{1}{\gamma(s,P_*)}\int_0^s\exp\left\{-\int_y^s\frac{\mu(z,P_*)}{\gamma(z,P_*)}\,\ud z\right\}g(y)\,\ud y\,\ud s<0,
\end{equation*}
holds true. In this case, we obtain for the growth bound of the semigroup $\omega_0$
\begin{equation*}
\omega_0=s(\mathcal{A+B+C+G}+\mathcal{H}_{\widetilde{\beta}})<0,
\end{equation*} 
see e.g. Theorem 1.15 in Chapter VI of \cite{NAG}, which implies that the equilibrium solution is linearly stable.
\end{example}

\section{Concluding remarks}

In this paper,  we analysed the asymptotic behaviour of a size-structured scramble competition model using linear semigroup methods. 
We are motivated by the modelling of structured macro-parasites in aquaculture, specifically the population dynamics of sea lice on Atlantic salmon populations. First we studied existence of equilibrium solutions 
of our model. In the case when the fertility function is separable, we easily established monotonicity conditions on the vital rates 
which guarantee the existence of a steady state (Proposition \ref{eq_separable}). In the general case we used positive 
perturbation arguments to establish criteria that guarantee the existence of at least one positive 
equilibrium solution. Next, we established conditions for the existence of a positive quasicontraction semigroup which governs the 
linearised problem. Then we established a further regularity property of the governing linear semigroup which in principle allows 
to study stability of equilibria via the point spectrum of its generator. In the special case of separable fertility function we 
explicitly deduced a characteristic function in equation \eqref{chareq} whose roots are the eigenvalues of the linearised operator. 
Then we formulated stability/instability results, where we used once more finite rank lower/upper bound estimates 
of the very general  recruitment term.  It would be also straightforward to formulate conditions which guarantee that the 
governing linear semigroup exhibits asynchronous exponential growth. However, this is not very interesting from the application point 
of view, since the linearised system is not necessarily a population equation anymore. 

Characterization of positivity using dispersivity resulted in much more relaxed conditions than those obtained in \cite{FH} for a more simple size-structured model with a single state at birth by characterizing positivity via the 
resolvent of the semigroup generator. This is probably due to the different recruitment terms in the two model equations.  
Positivity is often crucial for our stability studies, as was demonstrated in Section 3. Indeed, more relaxed positivity 
conditions result in the much wider applicability (i.e. for a larger set of vital rates) of our analytical stability results. 

Due to the fact that the positive cone of $L^1$ has an empty interior, characterizations of positivity such as the
positive minimum principle (see e.g. \cite{AGG}) do not apply. However, there is an alternative method, namely the
generalized Kato inequality (see e.g. \cite{AGG}). 
In our setting the abstract Kato-inequality reads
\begin{equation}\label{Kato}
S_u\,(\mathcal{A+B+C+D})u\le(\mathcal{A+B+C+D})|u|,
\end{equation}
for $u\in\text{Dom}(\mathcal{A+B+C+D})$, where $S_u$ is the signum operator, that is 
\begin{equation*}
S_u=\frac{u}{|u|}. 
\end{equation*}
Inequality \eqref{Kato} requires
\begin{align}
& S_u\int_0^m u(y)\left(\beta(s,y,P_*)+\int_0^m\beta(s,z,P_*)p_*(z)\,\ud z-\rho_*(s)\right)\,\ud y \nonumber \\
& \quad\quad\le\int_0^m|u(y)|\left(\beta(s,y,P_*)+\int_0^m\beta(s,z,P_*)p_*(z)\,\ud z-\rho_*(s)\right)\,\ud y, \quad s\in[0,m],\label{Kato2}
\end{align}  
which holds true for every $u\in\text{Dom}(\mathcal{A+B+C+D})$ indeed when condition \eqref{poscond1} is satisfied.

As we have seen previously in Section \ref{semigroup}, since the linearised system is not a population model anymore, 
the governing semigroup is not positive unless some additional condition is satisfied. However, it was proven in \cite{YK} that 
every quasicontraction semigroup on an $L^1$ space has a minimal dominating positive semigroup, called the modulus semigroup, which itself is 
quasicontractive. Hence, in principle, one can prove stability results even in the case of a non-positive governing semigroup, 
by perturbing the semigroup generator with a positive operator such that the perturbed generator does indeed generate a positive semigroup.

\section*{Acknowledgements}
JZF is thankful to the Centre de Recerca Mathem\`{a}tica and to the Department of Mathematics, 
Universitat Aut\`{o}noma de Barcelona for their hospitality while being a 
participant in the research programme ``Mathematical Biology: Modelling and Differential Equations'' during 01/2009-06/2009. 
PH thanks the University of Stirling for its hospitality. We also thank the Edinburgh Mathematical Society for financial support.


\end{document}